\documentclass[12pt]{article}
\usepackage{amsmath,amsfonts,amsthm,graphicx}

\theoremstyle{plain}
\newtheorem{thm}{\sc Theorem}[section]
\newtheorem{lem}[thm]{\sc Lemma}

\newtheorem{prop}[thm]{\sc Proposition}

\newtheorem{asser}[thm]{\sc Assertion}

\theoremstyle{remark}

\theoremstyle{definition}

\newcommand{\Si}{{\Sigma}}

\newcommand{\dc}{{\mathcal D}}

\newcommand{\lp}{{\mathcal P}}
\newcommand{\rth}{{\mathbb{R}^3}}
\newcommand{\W}{{\mathcal W}}

 \begin{title}
{The Topological Classification of Minimal Surfaces in
  $\rth$ }
\end{title}
\author{Charles Frohman and William H. Meeks III}

\begin{document}
\maketitle
\begin{abstract}We give a complete topological classification of properly
embedded minimal surfaces in Euclidian three-space
\end{abstract}
\thanks{This research was supported by NSF grant DMS -
 0104044 and NSF DMS 9803206.}
\section{Introduction.}

In 1980, Meeks and Yau \cite{my5} proved that properly embedded
minimal surfaces of finite topology in $\rth$ are unknotted in the
sense that any two such homeomorphic surfaces are properly ambiently isotopic.  
Later Frohman \cite{fr2} proved that any two triply periodic minimal surfaces are properly ambiently isotopic.  More recently
Frohman and Meeks \cite{fme1} proved that a properly embedded minimal
surface in $\rth$ with one end is a Heegaard surface
in $\rth$ and that Heegaard surfaces of $\rth$ with the same genus are unknotted; hence, properly
embedded minimal surfaces in $\rth$ with one end are unknotted even
when the genus is infinite.
These topological uniqueness theorems of Meeks and Yau, Frohman, and  Frohman
and Meeks are special cases of the following general classification
theorem which was conjectured in \cite{fme1} and which represents the final result for the topological classification problem.  The space of ends of a properly embedded minimal surface in $\rth$ has a natural linear ordering which is determined up to reversal and the middle ends in this ordering have a parity (even or odd) (see Section 2).

\begin{thm} (Topological Classification Theorem for
  Minimal Surfaces) \label{thm1.1} Two properly embedded minimal surfaces in $\rth$ are properly
  ambiently isotopic if and only if there exists a homeomorphism
  between the surfaces that preserves the ordering of their ends and preserves the
  parity of their middle 
  ends. 
\end{thm}

The constructive nature of our proof of the Topological Classification
Theorem provides an explicit description of any properly embedded minimal surface in terms of the
ordering of the ends, the parity of the middle ends, the genus of each
end - zero or infinite - and the genus of the surface.  This
topological description depends on several major advances in the classical theory of minimal surfaces.
First, associated to any properly embedded minimal surface $M$ with
more than one end is a unique plane passing through the origin called
the limit tangent plane at infinity of $M$ (see Section 2).
Furthermore, the ends of $M$ are geometrically ordered over
its limit tangent plane at infinity and this ordering is a
topological property of the ambient isotopy class of $M$ \cite{fme2}.  We call this result  the ``Ordering Theorem''.  Second, our proof of the 
classification theorem depends on the nonexistence of
middle limit ends for properly embedded minimal surfaces.   This result
follows immediately from the theorem of Collin, Kusner, Meeks and Rosenberg
\cite{ckmr1} that every middle
end of a properly embedded minimal surface in $\rth$ has quadratic
area growth.  Third, our
proof relies heavily on a topological description of the 
complements of $M$ in $\rth$; this topological description of the
complements was carried out by
the authors \cite{fme1} when $M$ has one end and by Freedman
\cite{fre1} in the general case.

Here is an outline of our proof of the classification theorem.  The
first step is
to construct a proper family $\mathcal{P}$ of topologically parallel
standardly embedded planes in $\rth$ such that the closed slabs and half spaces
determined by $\mathcal{P}$ each contains exactly one end of $M$ and
each plane in $\mathcal{P}$ intersects $M$ transversely in a simple
closed curve.  The
next step is to reduce
the global classification problem to a tractable
topological-combinatorial classification problem for ``Heegaard''
decompositions of closed slabs or half spaces in $\rth$.

\section{Preliminaries.}

Throughout this paper, all surfaces are embedded and proper.  We now
recall the definition of the limit tangent plane at infinity for a
properly embedded minimal surface $F\subset \rth$.  From the Weierstrass
representation for minimal surfaces one knows that the finite
collection of ends of a complete embedded
noncompact minimal surface $\Si$ of finite total curvature with compact boundary
are asymptotic to a finite collection of pairwise disjoint ends of planes and catenoids, each of which has
a well-defined unit normal at infinity.  It follows that 
the limiting normals to the ends of $\Si$ are parallel and one defines the limit tangent plane of
$\Si$ to be the plane passing through the origin and orthogonal to the
normals of $\Si$ at infinity.  Suppose that such a $\Si$ is contained in a complement of $F$.  One defines a limit tangent plane for $F$ to be
the limit tangent plane of $\Si$.  In \cite{chm3} it is shown
that if $F$ has at least two ends, then $F$ has a unique limit tangent
plane which we call {\it the limit tangent plane at infinity} for $F$.  We say that the limit tangent plane at infinity for
$F$ is horizontal if it is the $x_1 x_2$-plane.

The main result in \cite{fme2} is:

\begin{thm}(Ordering Theorem) \label{thm2.1} Suppose $F$ is a properly
  embedded minimal surface in $\rth$ with more than one end and with
  horizontal limit
  tangent plane at infinity.  Then the ends of $F$ have a natural
   linear ordering 
  by their ``relative heights'' over the $x_1 x_2$-plane.
  Furthermore, this ordering is topological in the sense that if
  $f$ is a diffeomorphism of $\rth$ such that $f(F)$ is a minimal
  surface with horizontal limit tangent plane at infinity, then the induced map on the spaces of ends preserves or reverses the orderings.
\end{thm}

Unless otherwise stated, we will assume that the limit tangent plane
at infinity of $F$ is horizontal and so $F$ is equipped with a
particular ordering on its set of ends $\mathcal{E}(F)$.
$\mathcal{E}(F)$ has a natural topology which makes it into a
compact Hausdorff space.  The limit points of $\mathcal{E}(F)$ are
called {\it limit ends} of $F$.  Since $\mathcal{E}(F)$ is compact and
the ordering on $\mathcal{E}(F)$ is linear, there exist unique maximal
and minimal elements of $\mathcal{E}(F)$ for this ordering.  The maximal element is called the {\it top} end of $F$.  The minimal element is called
the {\it bottom} end of $F$.  Otherwise the end is called a {\it middle} end of $F$.  

Actually for our purposes we will need to know how the ordering of
the middle ends $\mathcal{E}(F)$ is obtained.  This ordering is induced from a proper
family $\mathcal{S}$ of pairwise disjoint ends of horizontal planes and catenoids in
$\rth - F$ that separate the ends of $F$ in the following sense.
Given two distinct middle ends $e_1, e_2$ of $F$, for $r$ sufficiently large,
$e_1$ and $e_2$ have representatives in different components of
$\{(x_1, x_2, x_3) \in (\rth-\cup \mathcal{S})\mid x^2_1 + x^2_2 \geq r^2\}$.
Since the components of $\mathcal{S}$ can be taken to be disjoint
graphs over complements of round disks centered at the origin, they
are naturally ordered by their relative heights and hence induce an
ordering on $\mathcal{E}(F)$  \cite{fme2}.  

 In
\cite{ckmr1} it is shown that a limit end of $F$ must be a top or bottom end of the surface.  This means that each middle end $m \in \mathcal{E}(F)$
 can be represented by a proper subdomain $E_m \subset F$ which has
compact boundary and one end.  We now show how to assign a parity to
$m$.  First choose a vertical cylinder $C$ that
contains $\partial E_m$ in its interior.  Since $m$ is a middle end,
there exist components $K_+, K_-$ in $\mathcal{S}$ which are ends of
horizontal planes or catenoids in $\rth - F$ with $K_+$ above $E_m$
and $K_-$ below $E_m$.  By choosing the radius of $C$ large enough,
we may assume that $\partial K_+ \cup \partial K_-$ lies in the
interior of $C$.  Next consider a vertical line $L$ in $\rth - C$
which intersects $K_+$ and $K_-$, each in
a single point.  If $L$ is transverse to $E_m$, then $L\cap E_m$
is a finite set of fixed parity which we will call the parity of $E_m$.  The parity of $E_m$ only
depends on $m$, as it can be understood as the intersection number with
$\mathbb{Z}_2$-coefficients of the relative homology class of $L$,
intersected with the region between $K_+$ and $K_-$ and outside $C$,
with the homology class determined by the locally finite chain which
comes from the intersection of $E_m$ with
this same region.  If we let $A(R)$
denote the
 area of $E_m$ in the ball of radius $R$ centered at the
origin, then the results in \cite{ckmr1} imply that
$\lim_{R\rightarrow \infty} A(R)/\pi R^2$ is an integer with the same
parity as the end $m$. Thus, the parity of $m$ could also be defined geometrically
in terms of its area growth. This discussion proves the next
Proposition.

\begin{prop}\label{prop2.2} If $F$ is a properly embedded minimal
  surface in $\rth$, then each middle end of $F$ has parity.
\end{prop}

In \cite{fme1} Frohman and Meeks proved that the closures of the complements of a
minimal surface with one end in $\rth$ are handlebodies; i.e., they
are homeomorphic to the closed regular neighborhood of a properly
embedded connected $1$-complex in $\rth$.  Motivated by this
result and their ordering theorem, Michael Freedman \cite{fre1} proved the
following decomposition theorem for the closure of a complement of $F$ when $F$ has
possibly more than one end.

\begin{thm} \label{thm2.3}(Freedman)  Suppose $H$ is the closure of a complement of a
  properly embedded minimal surface in $\rth$.  Then there exists a
  proper collection $\mathcal{D}$ of pairwise disjoint minimal disks $(D_n,
  \partial D_n) \subset (H, \partial H), n \in \mathbb{N}$, such that the
  closed complements  of $\mathcal{D}$ in $H$ form a proper decomposition of
  $H$.  Furthermore, each component in this decomposition is a compact
  ball or is homeomorphic to $A \times [0, 1),$ where $A$ is an open annulus.
\end{thm}

\section {Construction of the family of planes $\lp$.}

\begin{lem} \label{lem3.1} Let $F$ be a properly embedded minimal
  surface in $\rth$ with one or two limit ends and horizontal limit tangent
  plane.  Suppose $H_1, H_2$ are the two closed
  complements of $F$ and $\dc_1$ and $\dc_2$ are the proper
  families of disks for $H_1, H_2$, respectively, whose existence is
  described in Freedman's Theorem.  Then there exist a
  properly embedded family $\lp$ of smooth planes transverse to $F$
  satisfying:

\begin{enumerate}
\item Each plane in $\lp$ has an end representative which is an end of
  a horizontal plane or catenoid which is disjoint
  from $F$;
\item In the slab $S$ between two successive planes in $\lp$, $F$ has
  only a finite number of ends;
\item Every middle end of $F$ has a representative in one of the just
  described slab
  regions $S$.
\end{enumerate}
\end{lem}

\proof The disks in $\dc_1$ can be chosen to be disks of least
  area in $H_1$ relative to their boundaries.  In fact the disks
  used by Freedman in the proof of his theorem have this property.  Assume that the disks in $\dc_2$ also have this least area
  property.  Suppose $W$ is a closed component of $H_1-\cup \dc_1$
  or $H_2-\cup \dc_2$ which is homeomorphic to $A \times [0, 1).$  Let
  $\gamma (W)$ be a smooth simple closed curve in $\partial W$ that
  generates the fundamental group of $W$. The curve of $\gamma (W)$ bounds a noncompact annulus in $\partial W$.
This annulus is a union of compact annuli $A_1 \subset A_2 \subset
\hdots A_n \subset \hdots$.  By \cite{my2} the boundary of $W$ is a
  good barrier for solving Plateau-type problems in $W$.  Let $\tilde{A}_n$ denote a least area annulus
in $W$ with the same boundary as $A_n$ which is embedded by \cite{my2}.
 The curve $\gamma (W)$ bounds
  a proper least area annulus $A(W)$ in $W$, where $A(W)$ is the limit
  of some subsequence of $\{\tilde{A}_n\}$; the existence of $A(W)$ depends on local curvature and local area estimates that we gave in a similar construction in \cite{fme1}.  By choosing $\gamma (W)$
  to intersect the interior of one of the disks in $\dc_1$ appearing
  in $\partial W$, the stable minimal annulus $A(W)$ intersects
  $\partial W$ only along $\partial A(W)$.  The stable minimal annulus
  $A(W)$ has finite total curvature \cite{fi1} and so is asymptotic to the end of a plane or catenoid in $\rth$.  By
  the maximum principle at infinity \cite{mr1}, the end of $A(W)$ is a positive
  distance from $\partial W$.  Hence, one can choose the
  representative end of a
  plane or catenoid to which $A(W)$ is asymptotic to lie in the
  interior of $W$.

Let $\mathcal{S}$ denote the collection of ends of planes and
catenoids defined above from all the nonsimply connected components $W$ of $H_i-\cup \dc_i$.  It
follows from the proof of the Ordering Theorem in \cite{fme2} that
$\mathcal{S}$ induces the ordering of
$\mathcal{E}(F)$.   Since the middle ends of $F$ are not limit ends, when $F$ has one limit end, then, after a possible
reflection of $F$ across the $x_1 x_2$-plane, we may assume that the limit
end of $F$ is its top end.  Thus, 
$\mathcal{S}$ will be naturally indexed by nonnegative integers $\mathbb{N}$ 
if $F$ has one limit end or by
$\mathbb{Z}$ if $F$ has two limit ends with the
ordering on the index sets 
$\mathbb{N}$ or $\mathbb{Z}$ coinciding with the natural ordering on $\mathcal{S}$ and the subset of  nonlimit ends in $\mathcal{E}(F)$.  

Suppose that $F$ has one limit end and let $S=\{ E_0, E_1, \hdots\}$.   Let $B_0$ be a ball of
radius $r_0$ centered at the origin with $\partial E_0 \subset B_0$ and such that $\partial B_0$
intersects $E_0$ transversely in a single simple closed curve $\gamma_0$.
The curve  $\gamma_0$ bounds a disk $D_0 \subset \partial B_0$.  Attach $D_0$ to
$E_0 - B_0$ to make a plane $P_0$.  Next let $B_1$ be a ball centered at the origin of
radius $r_1, r_1 \geq r_0 + 1$, such that $\partial E_1 \subset B_1$ and $\partial B_1$
intersects $E_1$ transversely in a single Jordan curve $\gamma_1$.
Let $D_1$ be the disk in $\partial B_1$ disjoint from $P_0$.  Let
$P_1$ be the plane obtained by attaching $D_1$ to $E_1 - B_1$.
Continuing in this manner we produce planes $P_n, n \in \mathbb{N}$,
that satisfy properties 1, 2, 3 in the Lemma.  If $F$ had two
limit ends instead of one limit end, then a simple modification of
this argument also would give a collection of planes $\mathcal{P}$
satisfying properties 1, 2, 3 in the Lemma.

\qed

\begin{prop} \label{prop3.2} There exists a collection of planes $\mathcal{P}$ satisfying the
  properties described in Lemma \ref{lem3.1} and such that each plane in $\mathcal{P}$
  intersects $F$ in a single simple closed curve.  Furthermore, in the
  slab between two successive planes in $\mathcal{P}$, $F$ has exactly one end.
\end{prop}  

\proof Suppose the limit tangent plane to $F$ is horizontal and
  that $\lp$ is finite.  Let $P_T$ and $P_B$ be the top and bottom
  planes in the ordering on $\lp$.  Since the inclusion of the
  fundamental group of $F$ into the fundamental group of either
  complement is surjective \cite{fme1}, the proof of Haken's lemma \cite{hak1} implies that
  $P_T$ can be moved by an ambient isotopy supported in a large ball
  so that the resulting plane $P^\prime_T$ intersects $F$ in a single
  simple closed curve.  Let $\tilde{P}_B$ be the image of $P_B$ under
  this ambient isotopy.  Consider the part $F_B$ of $F$ that lies in
  the half space below $P^\prime_T$ and note that the fundamental
  group of $F_B$ maps onto the fundamental group of each complement of
  $F_B$ in the half space.  The proof of Haken's lemma applied to $\tilde{P}_B$ in the
  half space produces an isotopic $P^\prime_B$ that intersects $F_B$
  in a simple closed curve.

Consider the slab bounded by $P^\prime_T$ and $P^\prime_B$.  The
following assertion implies that $\{P^\prime_T, P^\prime_B\}$ can be
expanded to a collection of planes $\mathcal{P}$ satisfying all of the
conditions of Proposition \ref{prop3.2}.

\begin{asser} \label{asser3.3} Suppose $S$ is a slab bounded by two
  planes in $\lp$ where $\lp$ satisfies Lemma \ref{lem3.1}.  Suppose each of
  these planes intersects $F$ in a simple closed curve.  Then there
  exists a finite collection of planes in $S$, each intersecting $F$
  in a simple closed curve, which separate $S$ into subslabs each of which contains a single  end of $F$.  Furthermore the addition of these planes to $\lp$
  gives a new collection satisfying Lemma \ref{lem3.1}.
\end{asser}

\proof Here is the idea of the proof of the assertion.  If $F$ has more than
one end in $S$, then there is a plane in $S$ topologically parallel to
the boundary planes of $S$ and which separates two ends
of $F \cap S$.  The proof of Haken's lemma then applies to give another such plane
with the same end
which intersects $F$ in a simple closed curve.  This new plane
separates $S$ into two slabs each containing fewer ends of $F$.  Since the number of
ends of $F\cap S$ is finite, the existence of the required
collection of planes follows by induction. 
\qed

Assume now that the number of planes in $\mathcal{P}$ satisfying Lemma \ref{lem3.1} is
infinite.  We first check that $\mathcal{P}$ can be refined to satisfy the
following additional property:~  If $W$ is a closed complement of either $H_1 -
\cup \mathcal{D}_1$ or $H_2 - \cup \mathcal{D}_2$, then $W$
intersects at most one plane in $\mathcal{P}$.  We will prove this in the case
that $F$ has one limit end.  The proof of the case where $F$ has two
limit ends is similar. 

Let $\mathcal{W}$ be the set of closures of the components of $H_1 - \cup \mathcal{D}_1$
and $H_2 - \cup \mathcal{D}_2$.  Given $W\in \mathcal{W}$, let
$\lp(W)$ be the collection of planes in $\lp$ that intersect $W$.  If
$W$ is a compact ball, then $\lp(W)$ is a finite set of planes since
$\lp$ is proper.  If $W$ is homeomorphic to $A \times [0,1)$, then
$\lp(W)$ is also finite.  To see this choose a plane $P \in \lp$ that
whose end lies above the end of $W$; the existence of such a plane is clear from the construction of $\lp$ in the previous Lemma.  Note that the closed half space above $P$
intersects $W$ in a compact subset.  Hence only a finite number of the
planes above $P$ can intersect $W$.  Since there are an infinite
number of planes in $\lp$ above $P$, there exists a plane $\tilde{P}$
above $P$ so that $\tilde{P}$ is disjoint from $W$ and any plane in
$\lp$ above $\tilde{P}$ is also disjoint from $W$.  Since there are
only a finite number of planes below $\tilde{P}$, only a finite number
of planes in $\lp$ can intersect $W$.  

We now refine $\lp$. First recall that the end of $P_0$ is contained in a
single component of $\mathcal{W}$. Hence, the plane $P_0$ intersects a finite number of
components in $\W$ and each of these components intersects a finite
collection of
planes in $\lp$ different from $P_0$.  Remove this collection from
$\lp$ and reindex to get a new collection $\lp = \{P_0, P_1, \cdots
\}$.  Note that $P_1$ does not intersect any component $W\in \W$ that
also intersects $P_0$.  Now remove from $\lp$ all the planes different
from $P_1$ that intersect some component $W\in \W$ that $P_1$
intersects.  Continuing inductively one eventually arrives at a
refinement of $\lp$ such that for each $W\in \W, \; \lp(W)$ has at most
one element.  This refinement of $\lp$ satisfies the conditions of
Lemma \ref{lem3.1} and so we may assume that $\lp(W)$ contains at most one plane
for every $P \in \lp$.

The next step in the proof is to modify each $P \in \lp$ so that the
resulting plane $P^\prime$ intersects $F$ in a simple closed curve.
We will do several modifications of $P$ to obtain $P^\prime$ and the
reader will notice that each modification yields a new plane that is a
subset of the union of the closed components of $\W$ that intersect
the original plane $P$.  This is important to make sure that 
further modifications can be carried out.

Suppose $P \in \lp$ and the end of $P$ is contained in $H_1$.  Let $\mathcal{A}_2$ be set of
components of $\W\cap H_2$ that are homeomorphic to $A\times[0,1)$.
For each $W \in \mathcal{A}_2$ let $T(W)$ be a properly embedded half
plane in
$W$, disjoint from $\cup \mathcal{D}_2$, such that the geodesic closure of $W
- T(W)$ is homeomorphic to a closed half space of $\rth$.  Assume that $P$ intersects transversely the half planes of the form $T(W)$
and the disks in $\mathcal{D}_2$.

We first modify $P$ so that there are no closed curve components in $P
\cap (\cup \mathcal{D}_2)$.  If $P \cap D, D \in \mathcal{D}_2$, has
a closed curve component, then there is an innermost one and it can be
removed by a disk replacement(see Figure A).  Since the end of $P$ is
contained in $H_1$, there are only a finite number of closed curve
components in $\cup \mathcal{D}_2$ and they can be removed by
successive innermost disk replacements.  In a similar way we can
remove the closed curve components in $P \cap ( \cup T(W))_{W \in
  \mathcal{A}_2}$.

We next remove compact arc intersections in $P \cap ( \cup \mathcal{D}_2)$ by
sliding $P$ over an innermost disk bounded by an innermost arc and
into $H_1$ (see Figure B).  In a similar way we can remove the finite
number of compact arc intersections of $P$ with $\cup T(W)_{W\in \mathcal{A}_2}$.
Notice that $P$ already intersects the region that we are pushing it
into.

After the disk replacements and slides described above, we may assume
that $P$ is disjoint from the disks in $\mathcal{D}_2$ and
half planes in $\mathcal{A}_2$.  Let $W \in \mathcal{W}$ be the component which contains the end of $P$ and let $P(*)$ be the component of $P
\cap W$ which contains the end of $P$.  Cut $H_2$ along the disks in
$\mathcal{D}_2$ and half planes in $\mathcal{A}_2$. Since every
closed component of the result is a compact ball or a closed
half space, the boundary curves of $P(*)$, considered as subsets of
the components, bound a collection of pairwise disjoint disks in
$H_2$.  The union of these disks with $P(*)$  is
a plane $P^\prime{^\prime}$ with  $P^\prime{^\prime} \cap W = P(*)$.  If $P(*)$ is an
annulus, then we are done.  Otherwise, since the fundamental group of
$W$ is $\mathbb{Z}$, the loop theorem implies that one can do surgery
in $W$
on $P(*) \subset P^\prime{^\prime}$ such that after the surgery, the component
with the end of $P^\prime{^\prime}$ has less boundary components.
After further surgeries in $W$ we obtain an annulus $P^\prime(*)$ with
the same end as $P(*)$ and with boundary curve being one of the
boundary curves of $P(*)$.  By our previous modifications
$\partial P^\prime(*)$ lies on the boundary of the closure of one of the components
of $H_2-\cup \mathcal{D}_2$ and bounds a disk $D$ in this
component.  We obtain the required modified plane $P^\prime =
P^\prime(*)\cup D$ which intersects $F$ in the curve $\partial P^\prime(*)$.

The above modification of a plane $P \in \lp$ can be carried out
independently since the modified plane is contained in the union of
the components of $\W$ that intersect $P$ and when  $P$ intersects
$W\in \mathcal{W}$, then no other plane in $\mathcal{P}$ intersects $W$.

Finally, applying the assertion at the beginning of the proof allows
one to subdivide the slabs between successive planes in $\lp$ so that
each slab contains at most one end of $F$.  This completes the
construction of $\lp$ and the proof of Proposition \ref{prop3.2}.
\qed

\section{The structure of a minimal surface in a slab.}

If $M$ is a 3-manifold and there is  a disjoint, properly
embedded system of disks $\mathcal{D}$ in $M$ so that the result of cutting
$M$ along $\mathcal{D}$ is a collection of balls, then $M$ is a 
{\em handlebody},
and $\partial M$ is the {\em preferred surface} of $M$. Alternatively, if $M$
is irreducible and there is a properly embedded $CW$-complex $\Gamma$ in $M$
so that $\Gamma$ is a strong deformation retract of $M$ and the
deformation $D\colon M \times [0,1] \rightarrow M$ is proper, then $M$ is a
handlebody. The second description is nice because we can perform handle-slides
and collapses on $\Gamma$ without changing the fact that it is a proper
deformation retraction of $M$.

We say $M$ is a {\em hollow handlebody} if there is a disjoint,  properly
embedded  system of disks $\mathcal{D}$ in $M$ so that the result of 
cutting $M$ along
$\mathcal{D}$ is homeomorphic to $\Sigma \times [0,1]$ for some surface
$\Sigma$ and  $\Sigma \times \{0\}$ lies completely in
$\partial M$. The surface $\partial M -\Sigma
\times \{0\}$ is the {\em preferred surface} of $M$.
If $H$ is irreducible and there is a subsurface $\Sigma$ of
$\partial H$ and a properly embedded $CW$-complex $\Gamma$ embedded in $H$ so
that  $\Sigma \cup \Gamma$ is a strong deformation retract of $H$ and the
deformation $D\colon H \times [0,1] \rightarrow H$ is proper, then $H$ is a hollow
handlebody.

There is yet another picture of handlebodies and hollow handlebodies that is
dual to the $CW$-complex $\Gamma$. Suppose $F$ is the preferred surface. There is a projection map 
$\pi\colon F \rightarrow \Gamma \cup \Sigma$.
The inverse image of $x \in \Gamma \cup \Sigma$ is
either a circle if $x$ is in the interior of an edge of $\Gamma$ or a
monovalent vertex of $\Gamma$, a theta
curve if $x$ is a trivalent vertex of $\Gamma$, or a point if $x \in \Sigma
-\Gamma$.

Choose a point in the interior of each edge of $\Gamma$. The inverse image of
this collection of points is a collection of circles. The circles decompose
$F$ into pairs of pants near trivalent vertices, annuli corresponding to
edges without trivalent vertices,  and a copy of $\Sigma$
with one disk removed for each monovalent vertex of $\Gamma$. We can
reconstruct a $CW$-complex that is isotopic to  $\Gamma$ from the system of 
circles.

\begin{figure}[htbp]
  \begin{center}
    \includegraphics{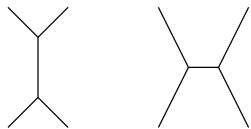}
    \caption{Whitehead move}
    \label{fig:white1}
  \end{center}
\end{figure}

\begin{figure}[htbp]
  \begin{center}
    \includegraphics{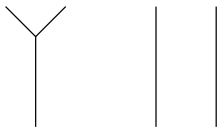}
    \caption{Half Whitehead move}
    \label{fig:white2}
  \end{center}
\end{figure}

Aside from isotopy there are two moves that we will be using on $\Gamma$.
They are both variants of the Whitehead move. We alter the graph according
to one of the two operations shown in Figure    \ref{fig:white1} and  Figure 
\ref{fig:white2}.

Dually the Whitehead move involves two pairs of pants meeting along simple
closed curve $\gamma$  which is the inverse image of a point in the interior
of the edge to be replaced. If $\gamma'$ is any simple closed curve lying on 
that union of pants that intersects $\gamma$ transversely in exactly two
points, 
and separates the boundary components of  the two pairs of pants into two 
sets of two, then we can
do the Whitehead move so that the two new pairs of pants meet along $\gamma'$.

\begin{figure}[htbp]
  \begin{center}
\includegraphics{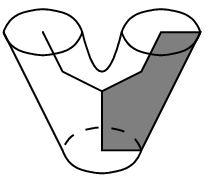} \hspace{1in}
\includegraphics{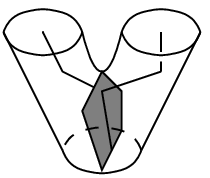}
\caption{Extending the disk $D$ to a singular surface.}
    \label{fig:huh}
  \end{center}
\end{figure}

Suppose that $H$ is a hollow handlebody or handlebody and $\delta$ is a
  simple closed curve in the preferred surface of $H$. We can extend $\delta$ to a
singular surface whose boundary lies in $\Gamma \cup \Sigma$. First isotope $\delta$ so
that with respect to the decomposition into annuli, pants and a punctured 
$\Sigma$, 
the part of $\delta$ that lies in each component is essential. 
There is a singular surface with boundary $\delta$ obtained
by adding ``fins'' going down to $\Gamma$ based on  the models
shown in Figure \ref{fig:huh}, along with fins  in the annuli and near $\Sigma$.

There are two kinds of essential arcs $k$ on a pair of pants, 
\begin{itemize}
\item {\bf good}: $\partial k$ runs between two distinct boundary components,

\item {\bf bad}: $\partial k$ joins a boundary component to itself.
\end{itemize}

\begin{lem}\label{propo}
  Suppose that $H$ is a handlebody or hollow handlebody and $\delta$ is a simple
  closed curve on the preferred surface of $ H$. Either $\delta$ bounds a disk in $H$ or there
is a graph $\Gamma$ so that $H$ is a regular neighborhood of 
$\Gamma \cup \Sigma$ such that $\delta$ has no bad arcs.
\end{lem}

\proof The argument will be by induction on a complexity for $\delta$.
Let $s$ be the number of bad arcs. Given a bad arc $k$, the  arcs (or arc)
of $\delta$  adjacent to $k$   lie in the same pair
of pants or in the punctured copy of $\Sigma$. If the two endpoints of the bad arc
coincide with the two endpoints of another bad arc, then let $d(k)=0$.
If both arcs lie in the punctured
copy of $\Sigma$, then let $d(k)=1$. If both arcs lie in the same pair of pants $P$, 
either the two arcs have their boundaries in the same components of $P$ or in
different components. If they have their endpoints in different components of
$P$, then $d(k)=1$. If not, then follow them into the next surface. If the next
surface is the punctured $\Sigma$, then $d(k)=2$, if the next surface is a pair of
pants and the next arcs are not parallel, then $d(k)=2$, otherwise follow them
into the next surface, and keep counting. Let $m=\min_{k \ \text{bad}}
{d(k)}$. The complexity of $\delta$ is the pair $(s,m)$.

\begin{figure}[htbp]
  \begin{center}
\begin{picture}(167,84)
    \includegraphics{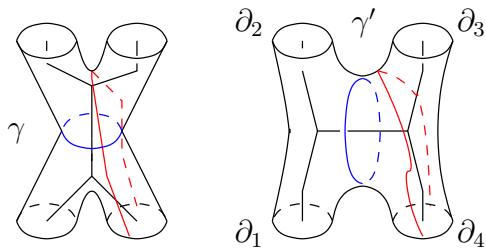} 
\put(0,0){$\partial_4$}\put(0,80){$\partial_3$}\put(-84,0){$\partial_1$}
\put(-84,80){$\partial_2$} \put(-40,80){$\gamma'$} \put(-170,40){$\gamma$}
\end{picture}
    \caption{Reducing $m$ when it is greater than $1$.}
    \label{fig:slide}
  \end{center}
\end{figure}

If $m>1$, then we do the Whitehead move as follows, see Figure
\ref{fig:slide}. Let $k$ be a bad arc with $d(k)=m$.
 Let $Q$ be the union of the pair of
pants containing $k$ and the pair of pants that contains the adjacent pair of 
arcs $k_1$ and $k_2$. Let $\gamma$ be the curve that that the two pairs of pants meet along.
Let $\partial_1,\partial_2,\partial_3,\partial_4$ be the boundary components
of $Q$ labeled so that $\partial_1$ and $\partial_4$ belong to one pair of
pants, $\partial_2,\partial_3$ belong to the other and both $k_1$ and $k_2$
have an endpoint in $\partial_4$.  
Let $a=k\cup k_1\cup k_2$. 
There is an arc $b$ of $\partial_4$ so that a push
off $\gamma'$ of $a\cup b$ lies in $Q$, misses $k$ and separates the boundary
components of $Q$ into two sets of two, say one set is $\partial_1$ and
$\partial_2$ and the other is $\partial_3$ and $\partial_4$.
Perform the Whitehead move so that
the pushoff  $\gamma'$ is the intersection of the two new pairs of pants.
Name the new pairs of pants $P_1$ and $P_2$. Notice that $a$ is a bad arc and
$d(a)=m-1$. 
To conclude that we have simplified the picture we need to see
that we have not increased the number of bad arcs. If $l$ is a bad arc in
$P_1\cup P_2$ and it has its endpoints in some $\partial_i$, then it contains
some bad arc of the original picture. If $l$ has its endpoints in $\gamma'$
and lies in $P_2$, as $\delta$ is embedded it is trapped in the annulus
between $\gamma'$ and $a$ and is hence inessential. If $l$ has its endpoints
in $\gamma'$ and is contained in $P_1$, once again the arc is trapped by $a$
and hence there must be two arcs in $P_2$ having one endpoint each in common
with $l$ and the other in $b$, but this means $l$ is contained inside a bad
arc from the original picture.
 Hence we did not
increase $s$. On the other hand we have decreased $m$ by $1$.

\begin{figure}[htbp]
  \begin{center}
    \includegraphics{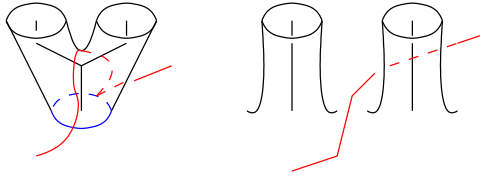}
    \caption{The endpoints of the bad arc are near $\Sigma$.}
    \label{fig:slide2}
  \end{center}
\end{figure}

If $m=1$, then there are two cases. The first is when next arcs lie in
the part of the surface parallel to $\Sigma$. In this case we do a half Whitehead
move and reduce the number of bad arcs, see  Figure \ref{fig:slide2}.

\begin{figure}[htbp]
  \begin{center}
    \includegraphics{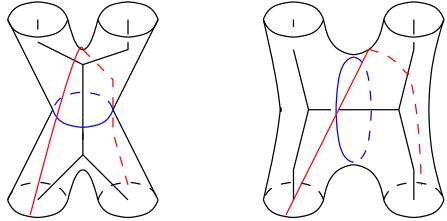}
    \caption{Reducing the number of bad arcs when $m=1$.}
    \label{fig:slide3}
  \end{center}
\end{figure}

The other case is when they are near another pair of pants. Once again, a
Whitehead move reduces the number of bad arcs, see Figure \ref{fig:slide3}.
Let $Q$ be the union of the two pairs of pants that contain $k$ and the
adjacent arcs $k_1$ and $k_2$, and let 
$\gamma$ be the circle that the pants intersect along.
Let $b$ be an arc in the pair of pants that does not contain $k$ and
which joins the
endpoints of $k$. Let $\gamma'$ be a pushoff of $b \cup k$ that
intersects $k$ in a 
single point and is disjoint from $k_1$ and $k_2$. Notice that the arc
$k_1\cup k \cup k_2$ gets separated into two good arcs by $\gamma'$. Hence if
we did not create any new bad arcs we have reduced the total number of bad
arcs. If a bad arc enters and leaves the new picture through a boundary
component of $Q$, then it is either contained in or contains a bad arc of the
old picture, hence we only need to worry about bad arcs with their endpoints
in $\gamma'$. Since $\delta$ is embedded such an arc misses $k_1 \cup k \cup
k_2$. The result of cutting $Q$ along the union of these arcs is a pair of
pants and $\gamma'$ gives rise to an arc of this pair of pants that has both
its endpoints in the same boundary component of the pair of pants. The only
proper arcs that intersect the arc corresponding to $\gamma'$ essentially in
two points must have both their endpoints in the same boundary component of
the pair of pants. This implies that such a bad arc is contained inside a bad
arc from the original picture.

\begin{figure}[htbp]
  \begin{center}
    \includegraphics{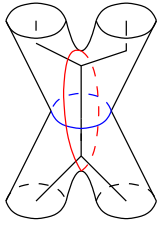}
    \caption{The interior disk.}
    \label{fig:slide4}
  \end{center}
\end{figure}

Finally when $m=0$ there are two arcs joined end to end, and the disk
inside the regular neighborhood of $\Gamma$ is readily visible, see  Figure
\ref{fig:slide4}.
\qed

A {\em  Heegaard splitting} with boundary of a $3$-manifold $M$ is a pair
$(H_1,H_2)$ where  $H_i$ is either a handlebody or hollow handlebody, and
$H_1 \cup H_2 =M$, and $H_1 \cap H_2$ consists of the preferred surface of
each $H_i$ . We call the preferred surface a {\em Heegaard surface} of $M$.

Suppose $S$ is a flat 3-manifold in $\mathbb{R}^3$ that is homeomorphic to
$\mathbb{R}^2 \times [0,1]$. Denote the components of $\partial S$ by
$\partial_0 S$ and $\partial_1 S$. Assume further that there are simple closed
curves $C_0 \subset \partial_0 S$ and $C_1 \subset \partial_1 S$ so that
$\partial_i S$ is a union of two stable minimal surfaces sharing $C_i$ as 
their joint
boundary.
As $\partial_i S$ is a plane,  one of these surfaces is a disk  $D_i$ and  the
other is a once punctured disk  $A_i$. Finally, assume that $F$ is a properly
embedded minimal surface in $S$ having one end and boundary $C_0 \cup C_1$.

\begin{prop}
The surface $F$ separates $S$ into two hollow handlebodies $H_1$ and $H_2$,
(or a handlebody and a hollow handlebody)
having $F$ as their preferred surfaces. That is, $F$ is a Heegaard
surface.

\end{prop}

\proof  We outline the idea for the sake of completeness.  
First consider the region $H_1$ and suppose that $\partial H_1$ has one end.  
In this case, by \cite{fme1} $H_1$ is a handlebody.  Assume now that 
$\partial H_1$ has two ends.  By Freedman's theorem applied to $H_1$, 
there exists a proper family of compressing disks $\dc_1$  which can be 
chosen to have their boundary components disjoint from $\partial S$.  After 
possibly restricting to a subcollection of $\dc_1$, the result of cutting 
$H_1$ along $\dc_1$ is connected and homeomorphic to $A \times [0,1)$.  But 
$A \times [0,1)$ is homeomorphic to $\Si \times [0,1]$ where $\Si$ is a 
proper once punctured disk on one of the boundary planes of $S$ with boundary one of 
the two boundary components of $F$.  In this case  $H_1$ is a hollow 
handlebody.  Similarly, if $\partial H_1$ has three ends, then one can 
choose the collection $\dc_1$ so that cutting $H_1$ along $\dc_1$ is 
homeomorphic to $\Si \times [0,1]$ where $\Si$ consists of the two once punctured disks 
in $\partial S$ bounded by $\partial F$.  Similarly, $H_2$ is either a 
handlebody or a hollow handlebody and so $F$ is a Heegaard surface in  $S$.
\qed

The proof of the topological classifiction theorem will require the 
examination of three kinds of surfaces with one end.

\begin{itemize}
\item {\bf Type 1.} The topology of $F \subset S$ is finite. This means that
  $F$ is homeomorphic to the result of removing a single point from a compact
  surface with two boundary components. In this case $F$ separates $S$ 
into two hollow handlebodies. One of
  the hollow handlebodies has boundary $D_0\cup F \cup A_1$ and the other has
  boundary $D_1\cup F \cup A_0$. Since $A_0$ and $A_1$ lie in different
  components of the complement of $F$, any arc joining $A_0$ to $A_1$ has
  $\mathbb{Z}_2$-intersection number $1$ with $F$. Hence the end is odd.

\item {\bf Type 2.} $F$ has infinite genus and any arc joining
  $A_0$ to $A_1$ has $\mathbb{Z}_2$-intersection number $1$ with $F$.  
Once again $F$ separates $S$ into two hollow handlebodies, one  with
  boundary  $D_0\cup F \cup A_1$ and  the other with boundary $D_1\cup F \cup
  A_0$ . This is an odd end.

\item {\bf Type 3.} $F$ has infinite genus and  any arc joining
  $A_0$ to $A_1$ has $\mathbb{Z}_2$-intersection number $0$ with $F$.
In this case $F$ separates $S$ into a handlebody with boundary $F \cup D_0
\cup D_1$ and a hollow handlebody with boundary $F \cup A_0 \cup A_1$. 
This end is even.
\end{itemize}

Our task is to show that in the first case, the surface is classified up
to topological equivalence  by its genus, and any two surfaces of the
second type (or third type) are topologically equivalent.
Let $D$ denote a topological disk, and let $A$ denote $S^1 \times [0,1)$.

\begin{thm} If $F$ and $F'$ are two minimal surfaces with one end  in $S$
  of finite type, the same genus and boundary consisting of circles $C_0$, $C_1$  and
  $C'_0$, $C'_1$ (respectively), then there
  is a homeomorphism $h\colon S \rightarrow S$ with $h(\partial_i S)=\partial_i S$ and $h(F)=F'$.
\end{thm}

\proof We will assume that we have chosen a homeomorphism between $S$ and
$\mathbb{R}^2 \times [0,1]$ and work in those coordinates. It is possible 
to find a large solid cylinder 
$D \times [0,1]$ whose boundary cylinder intersects $F$  in a single simple closed
curve in $\partial D \times [0,1]$ so that:
\begin{enumerate}
\item $S-D \times [0,1]$ is homeomorphic to $A \times [0,1]$;
\item The pair $(S-D \times [0,1], F-D\times [0,1])$ is topologically
  equivalent to the pair $(A \times [0,1],A \times \{1/2\})$.
\end{enumerate}

This follows quite easily from the fact that $F$ is  a Heegaard surface.
As $F$ has finite type, there is a {\em compact} $1$-dimensional CW-complex
$\Gamma$ so that $F$ is isotopic to the frontier of a regular neighborhood of
$\Gamma \cup \mathbb{R}^2 \times \{0\}$. Since $\Gamma$ is compact, its
projection to $\mathbb{R}^2$ is bounded. Hence there is a large $D$ in
$\mathbb{R}^2$ that contains its image. The set $D\times [0,1]$ satisfies
the conditions above. 
Similarly we could find $D' \times [0,1]$ having the same properties with
respect to $F'$.  

The existence of the $D$ above implies that we can compactify $S$, 
$F$ and $F'$ by adding a single circle at
infinity so that the compactification of $S$ is homeomorphic to the three-ball and the closures of $F$ and $F'$ are
Heegaard surfaces. The fact that $F$ and $F'$
complete to surfaces follows from the second property above. 
To see that  $F$ and $F'$ are Heegaard surfaces, note that the natural maps
on  fundamental groups induced by inclusion of the surfaces into their
complements are surjective. This implies that the compactified surfaces are
Heegaard splittings of the three-ball.  In \cite{fr3} it was proved
that such surfaces are classified up to homeomorphisms of the ball  by
their boundary and their genus. Hence if $F$ and $F'$ have the same genus, then we
can find a homeomorphism of the compactification of $S$ taking the
compactification of $F$ to the compactification of $F'$. By restricting the
homeomorphism we get a homeomorphism of $S$ having the desired properties.
\qed

Let $M$ be a manifold and suppose that $F$ is a Heegaard surface of
$M$ with compact boundary. We say that $F$ is {\it infinitely reducible} if
there is a properly embedded family of balls that are disjoint from
one another, so that each ball intersects $F$ in a surface of genus
greater than zero having a single boundary component, and so that
every end representative of $M$ has nonempty intersection with the
family of balls. It is a good exercise in the application of the
Reidemeister--Singer theorem to prove that any two infinitely
reducible Heegaard splittings of $M$ which
agree on the boundary of $M$
are topologically equivalent via a
homeomorphism of $M$ that is the identity on the boundary. This appeared
in \cite{fr4} and it can be seen from a proof analysis of \cite{fme1}.

Hence to show that up to topology there is only one surface in types $2$ and
$3$ it suffices to show that a minimal surface with one end of infinite
topology in a slab with boundary $C_0$, $C_1$ is infinitely reducible. For
this purpose we use a simple extension of a lemma from  \cite{fr2} to Heegaard
surfaces with boundary.

\begin{lem} Suppose that $F$ is the Heegaard surface of the irreducible 
manifold $M$, and there is a 1-dimensional $CW$-complex
  $\Gamma$ in $M$ and
  a subsurface $A$ of $\partial M$ so that $F$ is isotopic to a regular
  neighborhood of $\Gamma \cup A$. Suppose further that there is a ball $B$
  embedded in $M$ so that there is  a nontrivial cycle of $\Gamma$ 
 contained in the interior of $B$. Then $F$ is reducible.
\end{lem}

\proof Let $C$ be the nontrivial cycle of $\Gamma$ contained in the interior
of $B$. Notice that $F$ is a Heegaard surface for a splitting of the
complement of a regular neighborhood of $C$. Apply Haken's lemma to find a
sphere intersecting $F$ in a single circle. The sphere cuts off a subsurface 
of $F$
having genus greater than zero. Since the sphere bounds a ball in $M$, $F$ is
reducible. \qed

\begin{thm}
  If $F$ is a Heegaard surface of $S$ with one end, infinite genus and
  boundary consisting of two circles $C_i \subset \mathbb{R}^2\times \{i\}$, then the corresponding Heegaard splitting is infinitely
  reducible. \label{punch}
\end{thm}

\proof  of Theorem \ref{punch}.  Recall the coordinatization
$S=\mathbb{R}^2 \times [0,1]$. Let $\Gamma$ be a 1-dimensional
$CW$-complex so that
$F$ is the frontier of a regular neighborhood of $\Gamma$ and a subsurface of
$\partial (\mathbb{R}^2 \times [0,1])$, so that $\Gamma$ is of the form above.
By Proposition 2.2 of \cite{fme1}, there is an exhaustion of $S$ by compact submanifolds $K_i$ so that the 
part of
$F$ lying outside of each $K_i$ is a Heegaard surface for the complement of
$K_i$. 
For any $K_i$ there is $D_i \times [0,1]$ that contains $K_i$
and so that its frontier is transverse to $\Gamma$. Choose a half plane
$HP_i$ whose boundary consists of an arc in $\partial D_i \times [0,1]$ and
two rays, one each in $\mathbb{R}^2 \times \{0\}$ and $\mathbb{R}^2 \times
\{1\}$, that cuts the complement of $D_i \times [0,1]$ into a half space.
If there is a cycle of the graph $\Gamma$ in this half space, then there is a reducing ball outside $K_i$.

Since $F$ has infinite genus, there is a compressing disk $E$ for $F$ in the
complement of the regular neighborhood and that lies outside of $D_i \times
[0,1]$. By Lemma \ref{propo} there is either a disk that runs around $F$
without bad arcs or there is a disk inside the regular neighborhood
with boundary $\partial E$. In the second case the two disks form a sphere, which bounds a
ball in the complement of $D_i \times [0,1]$ containing a cycle of the graph.
In the second case, make $E$ transverse to $HP_i$. We can isotope $E$ (and the
graph $\Gamma$) so that there are no circles in $E \cap HP_i$. Let $k$ be an
outermost arc in $E\cap HP_i$. We will show that we can either alter the cycle which is the
boundary of $E$ so that it intersects $HP_i$ in fewer points or we can find a
nontrivial cycle of $\Gamma$ contained in the singular disk extending $E$.
In the case that we reduce the number of points we continue on. Either we find
a nontrivial cycle or we pull $E$ completely off of $HP_i$, in which case
there is a nontrivial cycle of $\Gamma$ disjoint from $HP_i$ in the desired
region. 

There are two cases. 
\begin{enumerate}
\item {\em The two endpoints of $k$ lie in the same boundary component of the
    same pair of pants.} 
 
The arc of the
  boundary of the disk extending E lying in $\Gamma$ defines a cycle. 
As $\partial E$ does not ever enter and leave a pair of pants
  through the same boundary component, this cycle is nontrivial. As the disk
  is outermost, there is a nontrivial cycle of $\Gamma$ in the result of
  cutting the complement of $D_i \times [0,1]$ along $HP_i$. As $HP_i$ is a half
  space it is easy to see there is a cycle of $\Gamma$ contained in a
  ball. Hence there is a trivial handle of $F$ lying outside $K_i$.

\item {\em The two endpoints of $k$ lie in distinct boundary components of
    pairs of pants.}

The first thing to notice is that the arc of the boundary of the disk
extending $E$ in 
$\Gamma$ is embedded. If not, then it would contain a cycle, and as the
disk is outermost that cycle would live in a ball. Let $l$ be the number of
pairs of pants that the arc passes through. If $l>1$, we reduce it via 
Whitehead moves on $\Gamma$ so as to not make any arcs of $\partial E$ bad. 
We can then use the outermost disk
as a guide to isotope $\Gamma$ so as to reduce the number of points of
intersection of that part of the graph in the boundary of the singular disk
which is the extension of $E$. It remains to show that if $l>1$, we can reduce
it. 

\end{enumerate}

After finitely many steps we have either found
a cycle in a ball or pulled the singular disk which is the boundary of $E$ off
of $HP_i$. In the case 2 above, because $E$ was a compressing disk, there is a
cycle of $\Gamma$ contained in the boundary of the singular disk, and it is
disjoint from $HP_i$ meaning we have a cycle in a ball and this ball lies outside of $K_i$ as desired.  This ball is contained in some $K_j, j>i$, and so we can reproduce this argument to find a cycle contained in a ball outside $K_j$.  It follows that the Heegaard splitting is infinitely reducible.
\qed

\proof Theorem  \ref{thm1.1}. Suppose that $F$ and $F'$ are two properly
embedded minimal surfaces and there exists a homeomorphism 
$h\colon F \rightarrow F'$ that preserves the ordering and parity of the ends.
By Proposition 3.2, we can find systems of planes that separate space into
slabs and the parts of $F$ and $F'$ lying in the respective slabs are Heegaard
surfaces. The parity and order preserving homeomorphism implies that there is
a correspondence between the slabs so that the parts of $F$ and $F'$ lying in
the corresponding slabs have the same parity. After shifting some handles
around so that finite genus surfaces have the same genus, we can then apply the
classification theorem for surfaces in a slab to build a homeomorphism
of the
space that takes $F$ to $F'$. \qed

\end{document}